\documentclass[12pt]{amsart}
\usepackage{amsmath,amsthm,epsfig,amssymb}
\title{Twin problems from the Monthly and the Stolz-Ces$\bf \grave{A}$ro Lemma }
\author{ Eugen J. Ionascu}

\date{\today}
\begin{document}
\def\sms{\small\scshape}
\baselineskip18pt
%%%%%%%%%%%%%%%%%%%%%%%%%%%%%%%%%%%%%%%%%%%%%%%%%%%%%%%%%%%%%%%%%%%%%%%%%%%%%%%%%%%%%%%%%%%%%%%%%%%%%
%%%%%%%%%%%%%%%%%%%%OUR DEFINITIONS%%%%%%%%%%%%%%%%%%%%%%%%%%%%%%%%%
%%%%%%%%%%%%%%%%%%%%%%%%%%%%%%%%%%%%%%%%%%%%%%%%%%%%%%%%%%%%%%%%%%%%%%%%%%%%%%%%%%%%%%%%%%%%%%%%%%%%%%%
\def\RR{{\rm I}\!{\rm R}}
\def\fp#1{(#1)}
\newtheorem{theorem}{\hspace{\parindent}
T{\scriptsize HEOREM}}[section]
\newtheorem{proposition}[theorem]
{\hspace{\parindent }P{\scriptsize ROPOSITION}}
\newtheorem{corollary}[theorem]
{\hspace{\parindent }C{\scriptsize OROLLARY}}
\newtheorem{lemma}[theorem]
{\hspace{\parindent }L{\scriptsize EMMA}}
\newtheorem{definition}[theorem]
{\hspace{\parindent }D{\scriptsize EFINITION}}
\newtheorem{problem}[theorem]
{\hspace{\parindent }P{\scriptsize ROBLEM}}
\newtheorem{conjecture}[theorem]
{\hspace{\parindent }C{\scriptsize ONJECTURE}}
\newtheorem{example}[theorem]
{\hspace{\parindent }E{\scriptsize XAMPLE}}
\newtheorem{remark}[theorem]
{\hspace{\parindent }R{\scriptsize EMARK}}
\renewcommand{\thetheorem}{\arabic{section}.\arabic{theorem}}
\newcommand{\lr}{L^2(\RR)}
\newcommand{\du}{\stackrel{.}{\bigcup}}
\renewcommand{\theenumi}{(\roman{enumi})}
\renewcommand{\labelenumi}{\theenumi}
\newcommand{\Q}{{\mathbb Q}}
\newcommand{\Z}{{\mathbb Z}}
\newcommand{\N}{{\mathbb N}}
\newcommand{\C}{{\mathbb C}}
\newcommand{\R}{{\mathbb R}}
\newcommand{\F}{{\mathbb F}}
\newcommand{\K}{{\mathbb K}}
\newcommand{\D}{{\mathbb D}}
\def\vp{\varepsilon}
\def\phi{\varphi}
\def\ra{\rightarrow}
\def\sd{\bigtriangledown}
\def\ac{\mathaccent94}
\def\wi{\sim}
\def\wt{\widetilde}
\def\bb#1{{\Bbb#1}}
\def\bs{\backslash}
\def\cal{\mathcal}
\def\ca#1{{\cal#1}}
\def\Bbb#1{\bf#1}
\def\blacksquare{{\ \vrule height7pt width7pt depth0pt}}
\def\bsq{\blacksquare}
\def\proof{\hspace{\parindent}{P{\scriptsize ROOF}}}
\def\pofthe{P{\scriptsize ROOF OF}
T{\scriptsize HEOREM}\  }
\def\pofle{\hspace{\parindent}P{\scriptsize ROOF OF}
L{\scriptsize EMMA}\  }
\def\pofcor{\hspace{\parindent}P{\scriptsize ROOF OF}
C{\scriptsize ROLLARY}\  }
\def\pofpro{\hspace{\parindent}P{\scriptsize ROOF OF}
P{\scriptsize ROPOSITION}\  }
\def\n{\noindent}
\def\iter#1{^{(#1)}}
\def\bh{\ca B(H)}
\def\ld{\overline}
\def\vpe{V_{\psi}^{\eta}}
\def\wu{\ca W(\ca U)}
\def \dper {\ca D\ca P}
\def\lir{L^{\infty}(\RR)} \def\ltr{L^2(\RR)} \def\ltrn{L^2(\RR^n)}
\def\lie{L^{\infty}(E)}\def\lte{L^2(E)}
\def\mpe{\ca M_{\psi,\eta}}
\def\wh{\widehat}
\def\eproof{$\hfill\bsq$\par}
\def\ws{\ca W\ca S}
\def\wi{\ca W\ca I}
\def\dst{\displaystyle}
\def\wsna{{\ca W}{\ca S}(n,A) }
\def\daw{{\ca W}_A}
\def\fwi{\ca W\ca I_1}
\def\swi{\ca W\ca I_2}
\def\ds{\displaystyle}
\def\du{\overset{\text {\bf .}}{\cup}}
\def\Du{\overset{\text {\bf .}}{\bigcup}}
\def\b{$\blacklozenge$}

%%%%%%%%%%%%%%%%%%%%%%%%%%%%%%%%%%%%%%%%%%%%%%%%%%%
%%%%%%%%%%%%%%%%%%%%%%%ABSTRACT%%%%%%%%%%%%%%%%%%%%%%%%%%%%%%%%%%%%
%\begin{abstract}
%In this note\footnotetext[1]{A version of this note appeared in
%...} we present solutions to two problems appeared which in the
%American Mathematical Monthly. Although it may seem that one
%problem is more general than  the other, actually the two problems
%cover different situations but both give sufficient conditions for
%a real valued function not to be periodic over the real line.
%However the two problems are related most intimately because they
%can be proved using essentially the same technique. In the end, we
%introduce a new problem which actually implies both.
%\end{abstract}
\maketitle

\n {\bf 1. INTRODUCTION.} In this note\footnote{A close version to
this note appeared in Crux Mathematicorum with Mathematical
Mayhem, vol. 34, Issue 7, (2008), pp. 424-429. We are posting this
here with the permission of the Canadian Mathematical Society.} we
present solutions to two problems which appeared in the American
Mathematical Monthly. Although it may appear that one problem is
more general than  the other, the two problems seem to cover
different situations but both give sufficient conditions for a
real valued function not to be periodic over the real line.
However the two problems are related most intimately because they
can be proved using essentially the same technique. In the end, we
introduce a new problem which actually implies both.  These
problems appeared in `04-`05 and their solutions in `06 and `07
(see \cite{pr1}, \cite{pr2}, \cite{pr1sol} and \cite{pr2sol}). The
proofs included here are based on a particular case of the
well-known Stolz-Ces$\grave{a}$ro Lemma and on the fact that every
continuous periodic function on $\mathbb{R}$ must be uniformly
continuous. The use of the latter idea is not new as it was used
in the published solutions of these problems. On the other hand,
the use of Stolz-Ces$\grave{a}$ro Lemma, is just another good
example where an old tool of analysis appears unexpectedly (see
\cite{batinetu},  \cite{ip}, \cite{italianu}, \cite{gnagy},
  and \cite{siretchi}).
L'Hospital's rule, which is very well known to calculus students
is its ``differentiable'' counterpart.

The version of Stolz-Ces$\grave{a}$ro Lemma we are going to employ
here is stated next.

\vspace{0.2in} {\bf Lemma 1.} {\it Let $\{a_n\}$ and $\{b_n\}$ be
two sequences such that $\{b_n\}$ is increasing and convergent to
infinity. If $\ds \lim_{n\to \infty} \frac{a_n}{b_n}=\infty$ then
$\ds \limsup_{n\to \infty} \frac{a_{n+1}-a_n}{b_{n+1}-b_n}=\infty$.}

\vspace{0.2in} We are going to include its classical idea of proof
for completeness. Let us assume to the contrary that $\gamma:=\ds
\limsup_{n\to \infty} \frac{a_{n+1}-a_n}{b_{n+1}-b_n}< \infty$.
Then for an arbitrary $\epsilon>0$ there exists $n_0\in \N$ such
that
$$\frac{a_{n+1}-a_n}{b_{n+1}-b_n}\le \gamma+\epsilon,$$
or

\begin{equation}\label{eq1}
a_{n+1}-a_n\le (b_{n+1}-b_n)(\gamma+\epsilon),
\end{equation}

\n for all $n\ge n_0$. Adding up inequalities as in (\ref{eq1}) for
$n=k...l$, $l>k\ge n_0$,  we obtain
$$a_{l+1}-a_k\le (b_{l+1}-b_k)(\gamma+\epsilon).$$
Eventually $b_{l+1}$ is going to be a positive number so we can
divide the last inequality by $b_{l+1}$ and then let $l\to
\infty$. Using the hypothesis we obtain $\infty \le
\gamma+\epsilon$ which is a contradiction.\eproof

Next, we are including the two original problems.

\par \vspace{0.1in}  {\bf Problem 11111 .} {\it  Let $f$
and $g$ be nonconstant, continuous periodic functions mapping $\R$
into $\R$. Is it possible that the function $h$ on $\R$ given by
$h(x)=f(xg(x))$ is periodic?} \vspace{0.1in}

\n The second problem seems to be more general but it is not clear
to us  at this point if this is indeed the case. We are going to
discuss the relationship between the two problems briefly but it
is not our purpose to get into the details of a thorough analysis.

\vspace{0.1in}  {\bf Problem 11174.} {\it  Let $f$ and $g$ be
nonconstant, continuous functions mapping $\R$ into $\R$ satisfying
the following conditions:\par \n 1. f is periodic.\par \n 2. There
is a sequence $\{x_n\}_{n\ge 1}$ such that $\ds \lim_{n\to
\infty}x_n=\infty$ and $\ds \lim_{n\to
\infty}\left|\frac{g(x_n)}{x_n}\right|=\infty$.\par \n 3. $f\circ g$
is not constant on $\R$.\par \n Determine whether $h=f\circ g$ can
be periodic.} \vspace{0.1in}

\n Both problems have a negative answer. The function
$h_1(x)=\sin(x\cos x)$ gives obvious choices for $f$ and $g$ that
satisfy the conditions in the first problem but it does not seem
to be an example (at least in an obvious way) good for the second
problem. On the other hand the function $h_2(x)=\sin(x^2)$ gives
rise to an $f$ and a $g$ that satisfy the conditions of the second
problem but it is hard to imagine that $h_2(x)={\hat f}(x{\hat
g}(x))$ for some $\hat f$ and $\hat g$ nonconstant, continuous
periodic functions. It is an interesting question whether or not,
for example, $h_1$ can be covered by Problem~11111.

The conditions in Problem~11174 can be weakened to obtain:

{\bf Theorem 1.} {\it Let $f$ and $g$ be nonconstant, continuous
functions mapping $\R$ into $\R$ and satisfying the following
conditions:\par \n (i) $f$ is periodic.\par \n (ii) there exist
sequences $\{x_n\}_{n\ge 1}$ and $\{y_n\}_{n\ge 1}$  such that
$$\ds \inf_{n} |x_n-y_n|>0\  \ and\ \
\lim_{n\to
\infty}\left|\frac{g(x_n)-g(y_n)}{x_n-y_n}\right|=\infty.$$\par \n
% 3. $f\circ g$ is not constant on $\R$.\par

\n Under these assumptions the function $h=f\circ g$ cannot be
periodic.}

\vspace{0.2in}

 \n {\bf 2. SOME FACTS FROM REAL ANALYSIS.}  Let us begin with this next fact about continuous
 functions on compact sets (Theorem 4.19 in \cite{rudin}).

{\bf Theorem 2.} {\it  Every continuous function on a compact set in
a metric space is uniformly continuous.}

We recall that a function on some domain $\ca D(f)\subset \Bbb R$ is
{\it uniformly continuous} if  for each $\epsilon >0$ there exists a
$\delta>0$ such that for every $x,y\in \ca D(f)$ for which
$|x-y|<\delta $  we have $|f(x)-f(y)|<\epsilon$.

An easy consequence of this theorem is:

{\bf Corollary 1.} {\it Every continuous periodic function on $\R$
is uniformly continuous.}

This can be seen by applying Theorem~2 to the restriction of a
periodic continuous function on $\R$ to the compact set $[0,2T]$
where $T>0$ is a period of $f$ and then taking the $\delta' $ given
by the theorem corresponding to an arbitrary $\epsilon$. This
$\delta'$ is good for the interval $[0,2T]$ as the theorem insures
but then $\delta:=min\{\delta',T\}$ is actually good for $\R$ as one
can easily verify.\par The idea of our proofs is to show that the
function $h$ is not uniformly continuous. As a result of Corollary~1
we see that $h$ cannot be periodic.

Let us see how Problem 11111 follows from Theorem~1. Since $g$ is
assumed to be continuous and periodic but not constant we can find
$a$ and $b$ such that $g(a)-g(b)\not =0$. Assume $T>0$ is a period
of $g$. Then we consider $x_n=a+nT$ and $y_n=b+nT$. Then
$|x_n-y_n|=|a-b|>0$ and

$$\lim_{n\to \infty} \frac{\left|x_ng(x_n)-y_ng(y_n)\right|}{|x_n-y_n|}=|a-b|^{-1}\lim_{n\to \infty}
\left|ag(a)-bg(b)+nT(g(a)-g(b)\right|= \infty,$$

\n which says that $x\overset{\hat{g}}{\to} xg(x)$ and $f$ satisfy
the conditions (i) and (ii) in Theorem~1  and so $f\circ \hat{g}=h$
is not periodic. So we have a solution for Problem 11111.

To show that Problem~11174 follows from Theorem~1 we need the weaker
version of the Stolz-Ces$\grave{a}$ro Lemma as stated in the
introduction as Lemma 1.

Now, let us assume the $f$, $g$ and $\{x_n\}$ satisfy the conditions
1-3 in Problem 11174. We can find a subsequence $\{x_{n_k}\}$ of
$\{x_n\}$, so that $x_{n_{k+1}}-x_{n_k}\ge 1$  for all $k$, and for
which either $\ds \lim_{k \to
\infty}\frac{g(x_{n_k})}{x_{n_k}}=\infty$ or $\ds \lim_{k \to
\infty}\frac{g(x_{n_k})}{x_{n_k}}=-\infty$. Without loss of
generality we may assume the first situation because the other case
is going to follow from this one by changing $g$ with $-g$ and $f$
with $x\overset{\hat f}{\to} f(-x)$ ($x\in \R$). By Lemma~1 we see
that
$$\ds \limsup_{k\to \infty} \frac{g(x_{n_{k+1}})-g(x_{n_k})}
{x_{n_{k+1}}-x_{n_k}}=\infty$$

\n which proves the existence of the two sequences in (ii) as in
Theorem~1. Hence, Theorem~1 can be applied to $f$ and $g$ and get
that $h=f\circ g$ is not periodic. This settles Problem 11174.

\vspace{0.3in}

 \n {\bf 3. PROOF OF THEOREM~1.} Let us start
with $f$ and $g$ satisfying (i) and (ii) of Theorem~1. Because $g$
is continuous and by property (ii) we see that the interval
$I_n:=g([x_n,y_n])$ (or $I_n:=g([y_n,x_n])$, for $n$ large enough,
must be an interval that has length greater than the period $T$ of
$f$. Hence the range of $f$ is the same as the range of $h=f\circ
g$. Since $f$ is assumed nonconstant then $h$ is nonconstant.
Therefore we can choose $\alpha$ and $\beta$ such that
$f(g(\alpha))\not=f(g(\beta))$ and then we let
$\epsilon_0=|f(g(\alpha))-f(g(\beta))|>0$. As we said in the
introduction the key idea is to prove that $h$ is not uniformly
continuous. More precisely, we want to show that the definition of
uniform continuity is not satisfied for this $\epsilon_0$.

%According to Lemma~\ref{scl} and the comment we just made, we obtain
%the existence of two sequences $\{y_n\}$ and $\{z_n\}$ convergent to
%infinity, both increasing,  and such that

%\begin{equation}\label{ec}
%z_n-y_n\to \infty \ \text{and}\ \ \lim_{n\to \infty}
%\frac{g(z_n)-g(y_n)}{z_n-y_n}=\infty.
%\end{equation}
%This is saying roughly speaking that the size of the image of the
%interval $[y_n,z_n]$ through $g$ relative to the size of $[y_n,z_n]$
%is growing to infinity.
%We are going to denote by $T$ ($T>0$) a period of $f$.

We fix $n\in \N$ large enough to insure that $|I_n|>2T$ and denote
by $\sharp(g(\alpha))$ the number of integer values of $k$ for which
$g(\alpha)+kT$ is in $I_n$. Then, it is easy to see that
$$\ds \sharp(g(\alpha))> \frac{|g(x_n)-g(y_n)|}{T} -1>1.$$

Similarly we denote by $\sharp(g(\beta))$, the number of  integers
$k$ for which $g(\beta)+kT$ is in $I_n$. Again, we have
$\sharp(g(\beta))>\frac{|g(x_n)-g(y_n)|}{T}-1>1$.

It is clear that the values $g(\alpha)+kT$ ($k\in\Z$) interlace with
those of $g(\beta )+kT$ ($k\in\Z$). Using again the fact that $g$ is
continuous, by repeated application of the Intermediate Value
Theorem we can find two finite sequences ${u_k}$ and ${v_k}$ in the
interval $[x_n,y_n]$ (or $[y_n,x_n]$) both increasing and
interlacing such that $g(u_k)=g(\alpha)+l_kT$ and
$g(b_k)=g(\beta)+s_kT$ with $l_k,s_k\in \Z$. The number of the
intervals of the form $[u_k,v_k)$ (or $[v_k,u_k)$, $[v_k,u_{k+1})$,
etc.) is at least
$$M:=min[2(\sharp(g(\alpha)-1),2(\sharp(g(\beta)-1)]\ge 2.$$

These intervals form a partition of a subinterval of
$J_n:=[x_n,y_n]$ (or $J_n:=[y_n,x_n]$) of length $|x_n-y_n|$. It
follows that at least one of these intervals has to have length less
than or equal to $\frac{|x_n-y_n|}{M}$.

We denote such an interval by $[\zeta_n,\eta_n]$ and notice that

\begin{equation}\label{ess}
\ds |\zeta_n-\eta_n|\le \frac{|x_n-y_n|}{M}<
\frac{|x_n-y_n|}{2\frac{|g(x_n)-g(y_n)|}{T}-4}=
\frac{1}{\frac{2}{T}\frac{|g(x_n)-g(y_n)|}{|x_n-y_n|}-\frac{4}{|x_n-y_n|}}\to
0\ as\ n\to \infty,
\end{equation}

\n and $|f(g(\zeta_n))-f(g(\eta_n)|=\epsilon_0$. For an arbitrary
but fixed $\delta>0$, we choose $n$ even bigger so that
$|\zeta_n-\eta_n|<\delta$. This can be done because of (\ref{ess}).
For such an $n$ we still have $|h(\zeta_n)-h(\eta_n)|\ge \epsilon_0$
which proves that $h$ is not uniformly continuous.

In the end we would like to leave the reader with a natural
question: can Theorem~1 be generalized to almost periodic functions?
There are various concepts of almost periodicity but we are going to
include here as an example only Bohr's definition:

{\it A continuous real valued function $\digamma$ defined on $\R$ is
said to be {\bf almost periodic} if for each $\epsilon>0$ there
exists an $L>0$ such that every interval of length $L$ contains an
$\epsilon$-period, i.e. a number $T$ such that
$|\digamma(x+T)-\digamma(x)|<\epsilon$ for all $x\in \R$. }

\n What we find encouraging when it comes to new developments,
related to the above questions, is the fact that every almost
periodic function is also uniformly continuous (\cite{cor}).

\n {\bf Acknowledgements:} We thank professor Albert VanCleave who
gave us helpful suggestions after reading an earlier version of
this note.

\small Department of Mathematics, Columbus State University\\
\small Columbus, GA 31907; \texttt{ionascu\_eugen@colstate.edu}\\
\small Honorific Member of the Romanian Institute of Mathematics ``Simion Stoilow''\\


\begin{thebibliography}{99}

\bibitem{aops} Art of problem Solving, {\em
http://www.artofproblemsolving.com/}




 \bibitem{batinetu} D. M. Batinetu-Giurgiu, Siruri (Sequences), Editura Albatros, Bucuresti 1979.


\bibitem{cor} C. Corduneanu,   Almost Periodic Functions. New York: Wiley
Interscience, 1961.
\bibitem{pr1} P.P. Dalyay, Problem 11111, {\em Amer. Math. Monthly}, Vol.
111, {\bf no. 9}, 2004, p. 822.
\bibitem{pr2} P.P. Dalyay, Problem 11174, {\em Amer. Math. Monthly}, Vol.
112, {\bf no. 8}, 2005, p. 749.

\bibitem{ip} E. J. Ionascu and  P. Stanica, {\em Effective Asymptotics for Some Nonlinear
Recurrences and Almost Doubly-Exponential Sequences}, Acta Math.
Univ. Comenian. (N.S.) {\bf no. 73} (2004), no. 1, 75--87.

\bibitem{pr1sol} GCHQ Problem Solving Group, {\em Amer. Math. Monthly}, Vol.
113, {\bf no. 5}, 2006, p. 467.

\bibitem{italianu} A. Mannino, {\em Some classic Stolz-Cesàro theorems.
(Italian)}, Atti Accad. Sci. Lett. Arti Palermo Parte I (4) 37
(1977/78).
\bibitem{gnagy} G. Nagy, {\em The Stolz-Cesaro Theorem} http://www.math.ksu.edu/~nagy/snippets/stolz-cesaro.pdf
\bibitem{pr2sol} National Security Agency Problems Group, {\em Amer. Math. Monthly}, Vol.
114, {\bf no. 9}, 2007, p. 836.
\bibitem{kr} K. A. Ross, {\em Elementary Analysis: The theory of
calculus}, Springer-Verlag, 1980
\bibitem{rudin} W. Rudin, {\em Principles of Mathematical Analysis}, Third Edition, McGraw-Hill, Inc. 1964.
\bibitem{siretchi} Gh. Sire\c tchi, {\em The Toeplitz theorem and some of its
consequences (Romanian)}, Gaz. Mat. (Bucharest) 90 (1985), no. 3,
65--70.
\end{thebibliography}
\end{document}